\documentclass[12pt, twoside, leqno]{article}

\usepackage[T1]{fontenc}
\usepackage{multirow}
\usepackage{amsmath,amsthm}
\usepackage{amssymb}
\usepackage[utf8]{inputenc}
\usepackage{color}
\usepackage{booktabs}

\usepackage{enumerate}

\usepackage{graphicx}

\newtheorem{thm}{Theorem}[section]

\newtheorem{lem}[thm]{Lemma}

\theoremstyle{definition}
\newtheorem{defin}[thm]{Definition}

\theoremstyle{preuve}

\numberwithin{equation}{section}

\begin{document}


\baselineskip=17pt


\title{Number of points on a family of curves over a finite field}

\author{Thiéyacine Top${}^{*}$}
\date{}
\maketitle
\renewcommand{\thefootnote}{}
\footnote{$^{*}$ Université Blaise Pascal, Clermont-Ferrand, France
(E-mail: thieyacine.top@math.univ-bpclermont.fr)}
\setcounter{footnote}{0}


\begin{abstract}

In this paper we study a family of curves obtained by fibre products of hyperelliptic curves. We then exploit this family
to construct examples of curves of  given genus $g$  over a finite field $\mathbb{F}_{q}$ with many rational points.
The results obtained improve the known bounds for a few pairs $(q,g)$.\\

\textbf{Key words.} Curves over finite fields, jacobian varieties.\\
\textbf{AMS subject classification.} 11G20, 14G15, 14H45.
\end{abstract}

\section{Introduction}
Let  $\mathcal{C}$ be a (projective, non-singular and geometrically irreducible) curve of genus $g$ defined over the finite field $\mathbb{F}_{q}$
with $q$ elements. We denote by $N_{q}(g)$  the maximum number of rational points on a curve of genus $g$ over $\mathbb{F}_{q}$, namely
$$ N_{q}(g) = \max \left\{ | \mathcal{C}(\mathbb{F}_{q})| : \mathcal{C}~~ \text{is a curve over } \mathbb{F}_{q} ~~\text{of genus}~~ g \right \}.$$
In the last years, due mainly to applications in Coding Theory and Cryptography (see e.g.[8]) , there has been considerable interest in computing $N_{q}(g)$.
It is a classical result that $N_{q}(0) = q + 1$. Deuring and Waterhouse[10], and Serre [6] computed $N_{q}(1)$  and $N_{q}(2)$ respectively. Serre also computed $N_{q}(3)$ for $q < 25$  and Top [7] extended these computations to $q < 100$.

For $g \geq 3 $ no such general formula is known. However, Serre [6] building on earlier results by Hasse and Weil proved the following
upper-bound:
$$  N_{q}(g) \leq q + 1 + g \lfloor 2\sqrt{q} \rfloor ,$$ where $\lfloor x \rfloor$ denotes the floor part of a real number $x$.
This bound is now called Hasse-Weil-Serre bound.\\
Some partial results for some specific pairs $(q,g)$ are recorded on the website
\texttt{http:// www.manypoints.org} in the form $N_{q}(g) \in [c_{1}, c_{2}]$ or $[\ldots,c_{2}]$, where
\begin{itemize}
  \item $c_{2}$ is the bound given by Hasse-Weil-Serre or by  "explicit formulas", or by more intricate arguments.
  \item $c_{1}$ is the bound given by the existence of a curve  $\mathcal{C}$ with $ |\mathcal{C}(\mathbb{F}_{q})|= c_{1}$,
  \item \ldots~~ when $c_{1} \leq c_{2}/ \sqrt{2}$.
\end{itemize}
Our work consists in  studying   geometric families of curves and studying among them those which have many points over $\mathbb{F}_{q}$.
The curves we are considering are certain fibre products of hyperelliptic curves, we prove a formula for the genus and a formula for the number
of rational points of such curves that depend are the polynomial $f_{i}$'s and the hyperelliptic they define.
In particular by building  curves of given genus explicitly, we have the following result:\\
Let $k$ be an integer $\geq 1$ and $q \geq 2$ be a prime power. Let $I$ be a non-empty subset $\{1,2,\ldots, k\}$ and denote by $\mathcal{I}$
the set of all non-empty subsets of $I$. Let $f_{1}, f_{2}, \ldots , f_{k}$ be polynomials of respective degrees $d_{i}$'s, with coefficients in the finite fields $\mathbb{F}_{q}$
such that the polynomial $f_{1}\times f_{2}\times \ldots \times f_{k}$ is separable. We define the polynomial
$$f_{I}(x) = \prod_{i \in I} f_{i}(x) $$ and denote by $\mathcal{C}_{I}$ the hyperelliptic curve of equation $y^{2} = f_{I}(x)$. Let $A_{I} = q + 1 - | \mathcal{C}_{I}(\mathbb{F}_{q})|$. Consider now the fibre product $\mathcal{C}$ of which an  open affine subset is given by
$$\left \{ \begin{array}{ccc}
  y_{1}^{2} & = & f_{1}(x)  \\
   & \vdots &  \\
  y_{k}^{2} & = & f_{k}(x)
\end{array} \right.,$$
Let $N = |\mathcal{C}(\mathbb{F}_{q})|$ where $|\mathcal{C}(\mathbb{F}_{q})|$ is the number of $\mathbb{F}_{q}$-rational points on $\mathcal{C}$.
\begin{thm} In this notation:
\begin{itemize}
  \item The genus $g$ of $\mathcal{C}$ is givin by: $$g = 2^{k-2}(d_{1} + \ldots + d_{k} -4) +1 + \delta_{k}, $$
with $\delta_{k} = 2^{k-2}$ (resp. $\delta_{k} = 0$) if one of the $d_{i}$'s is odd (resp. all of the $d_{i}$'s are even). There are $2^{k-1}$
 (resp. $2^{k}$) points at infinity  if one of the $d_{i}$'s is odd (resp. all of the $d_{i}$'s are even).
  \item The number of $\mathbb{F}_{q}$-rational points of $\mathcal{C}$ is given by:
  $$N = q +1 - \sum_{I \in \mathcal{I}}A_{I}.$$
  \end{itemize}
\end{thm}
The proof of the above theorem follows from the lemma 2.1. and the lemma 3.3. for the computation of $g$ and $|\mathcal{C}(\mathbb{F}_{q})|$ respectively.
In a second part of the paper, we use Theorem1.1 on specific examples in order to improve a current lower bounds for $N_{q}(g)$ for some
values of $(q,g)$. The results we obtain are listed in Table 1.\\
\begin{table}[h]
   \begin{tabular}{|c|c|c|c||c|c|c|c|}
         \hline
         \multirow{2}{0.5cm}{$g$} & \multirow{2}{0.5cm}{$q$} &
         \multicolumn{2}{c|}{$N_{q}(g)$} &  \multirow{2}{0.5cm}{$g$} & \multirow{2}{0.5cm}{$q$} &
         \multicolumn{2}{c|}{$N_{q}(g)$}\\
         \cline{3-4} \cline{7-8}
         & & New enter & Old &&&New enter & Old\\
         \hline
         5 & $17$ & 48 & $[\ldots ,53]$  & 5 & $73$ & 148  & $[\ldots ,156]$ \\
           & $19$ & 52& $[\ldots ,60]$ & & $79$ & 156 & $[\ldots ,165]$ \\
          & $23$ & 62 & $[\ldots ,67]$ & & $83$ & 162  & $[\ldots ,172]$ \\
          & $29$ & 72  & $[\ldots ,80]$ &  & $89$ & 168 & $[136 ,180]$  \\
          & $31$ & 76& $[\ldots ,84]$ & & $97$ & 180 &  $[\ldots ,193]$ \\
          & $37$ & 88  & $[\ldots ,96]$  &  & $5^2$ & 64 & $[\ldots ,72]$ \\
          & $41$ & 94  & $[\ldots ,102]$  & &$13^2$ & 295  & $[232 ,300]$ \\
          & $43$ & 100 & $[\ldots ,106]$   & &$17^2$ & 454  & $[376 ,460]$ \\\cline{5-8}
          & $47$ & 102  & $[\ldots ,113]$  & 6 & $23$ & 66 & $[60 ,78]$ \\
          & $53$ & 120 & $[\ldots ,124]$   &  & $31$ & 84 & $[80 ,92]$\\
          & $59$ & 124 &  $[\ldots ,133]$ & & $41$ & 104& $[102 ,114]$\\
          & $61$ & 126 & $[\ldots ,137]$  &  & $59$ & 134 & $[132 ,150]$ \\\cline{5-8}
          & $67$ & 136  & $[\ldots ,148]$  & $7$ & $29$ & $80$ & $[72 ,100]$ \\ \cline{5-8}
          & $71$ & 144  & $[\ldots ,152]$ & $8$ &11& $46$ &   $[42 ,55]$  \\\cline{5-8}
         \hline
      \end{tabular}
      \caption{The specified interval (Old) is that given by  the website
\texttt{http:// www.manypoints.org} as of Avril 2016}
      \end{table}

\section{Geometry of curves}
Let $k$ be a positive integer. Consider now the fibre product $\mathcal{C} = \mathcal{C}_{f_{1} \ldots f_{k}}$ of which an  open affine subset is given by
$$\left \{ \begin{array}{ccc}
  y_{1}^{2} & = & f_{1}(x)  \\
   & \vdots &  \\
  y_{k}^{2} & = & f_{k}(x)
\end{array} \right.,$$

where the $f_{i}$'s are coprime polynomials over $\mathbb{F}_{q}$ of respective degrees $d_{i}$'s.
\begin{lem}Suppose the polynomial $f(x) := \prod_{i}f_{i}(x)$  is separable, then the affine curve is smooth. The genus of the complete curve is $$ g = 2^{k-2}(d_{1} + \ldots + d_{k} -4) +1 + \delta_{k},$$
with $\delta_{k} = 2^{k-2}$ (resp. $\delta_{k} = 0$) if one of the $d_{i}$'s is odd (resp. all of the $d_{i}$'s are even). There are $2^{k-1}$
 (resp. $2^{k}$) points at infinity  if one of the $d_{i}$'s is odd (resp. all of the $d_{i}$'s are even).
 \end{lem}
\textbf{Proof.} Smoothness follows from the jacobian criterion applied to the matrix $$\left( \begin{array}{ccccc}
             f'_{1}(x) & 2y_{1} & 0 &  &  \\
             f'_{2}(x) & 0 &  2y_{2}& 0 &  \\
              &   &  & \ldots &  \\
              f'_{k}(x)&  &  & 0 & 2y_{k}
           \end{array}
\right ).$$

 If none of the $y_{i}$'s is zero, there exists a minor equal to $$2^{k-1}y_{1} \ldots y_{i-1}f'_{i}(x)y_{i+1}\ldots y_{k} \neq 0; $$
If for some $i$, $y_{i} = 0$, then $f_{i}(x) = 0$ and therefore $f'_{i}(x)\neq 0$ and, for $j \neq i$, we have $f_{j}(x)\neq 0$ and there exists a minor equal to $$2^{k-1}y_{1} \ldots y_{i-1}f'_{i}(x)y_{i+1}\ldots y_{k} \neq 0 .$$
The group $G = \left\{\pm 1\right \}^{k} \cong (\mathbb{Z}/2\mathbb{Z})^{k}$ acts in an obvious way on $\mathcal{C}$ by $$[\varepsilon](x,y_{1}, \ldots , y_{k})= (x,\varepsilon_{1}y_{1}, \ldots , \varepsilon_{k}y_{k})$$  and is the Galois group of the covering $$\phi : \mathcal{C} \rightarrow \mathbb{P}^{1}$$ given by $$(x,y_{1},\ldots, y_{k})\mapsto x.$$ The group  $G$ acts transitively on the set $\mathcal{C}_{\infty}$, of points at infinity; the inertia group is cyclic therefore is trivial or reduced to $\mathbb{Z}/2\mathbb{Z}$.\\
If one of the  $d_{i}$'s is odd (resp.  all $d_{i}$'s are even), there is ramification above $\infty \in \mathbb{P}^{1}$ hence we obtain $|C_{\infty}|= 2^{k-1}$ (resp. $= 2^{k}$ ).
The  Riemann-Hurwitz formula applied to the morphism $\phi$ gives the genus, observing that $\phi$ is branched at every point $x = \alpha_{i}$
with $f_{i}(\alpha_{i}) = 0$ , i.e $$(x,y_{1},\ldots,y_{k}) = (\alpha_{i},\pm \sqrt{f_{1}(\alpha_{i})}, \ldots,0, \ldots ,\pm \sqrt{f_{k}(\alpha_{i})}),$$ and possibly above $\infty$.
\section{Computation  the number of points}
We keep the same notation as before.
\begin{defin}Let $I$ be a non-empty subset of $\{1, \ldots, k\}$, we define:
\begin{enumerate}
  \item the polynomial $f_{I}(x) = \prod_{i \in I}f_{i}(x)$ and we denote by $d_{I}$ its degree;
  \item the smooth projective curve $\mathcal{C}_{I}$  whose affine model is given by $v^{2} = f_{I}(u)$;
  \item we denote by $g_{I}$ the genus of $\mathcal{C}_{I}$;
  \item the morphism $\phi_{I}: \mathcal{C} \rightarrow C_{I}$ given by  $\phi_{I}(x,y_{1}, \ldots , y_{k}) = (x,y_{I})$ (where $y_{I} = \prod_{i \in I}y_{i})$.
\end{enumerate}
We denote by $\mathcal{I}$ the set of non-empty subsets of $\{1, \ldots, k\}$. \end{defin}

\begin{lem} The morphism
$$
\begin{array}{rcl}
  \Psi: \prod_{I \in \mathcal{I}} \mathrm{Jac}(\mathcal{C}_{I})& \longrightarrow  & \mathrm{Jac}(\mathcal{C}) \\
  (\mathcal{C}_{I})_{I \in \mathcal{I}} & \longmapsto  &\sum_{I \in \mathcal{I}}\phi_{I}^{\ast}(\mathcal{C}_{I }) .
\end{array}
 $$ is a separable isogeny.  \end{lem}
\textbf{Proof.} We first verify that $\sum_{I \in \mathcal{I}}g_{I} = g$; indeed, if we put $d_{I} = \sum_{i \in I} d_{i} = \mathrm{deg} f_{I}$
then $g_{I} = \left \lfloor\frac{d_{I}-1}{2} \right\rfloor $, which we may write $\frac{d_{I}-1-\varepsilon_{I}}{2}$
with $\varepsilon_{I} = 0$ or 1. \\ So
$$\sum_{I \in \mathcal{I}}g_{I} = \sum_{I \in \mathcal{I}}\frac{d_{I}-1-\varepsilon_{I}}{2}= \frac{1}{2} \left( \sum_{i = 1}^{k}
d_{i}N_{i} - 2^{k}+1 - \sum_{I}\varepsilon_{I}\right)$$
 where $N_{i}$ is the number of  $I$'s containing $i$, i.e $N_{i} = 2^{k-1}$.\\ If all $d_{i}$ are even, all the $\varepsilon_{I}$ = 1 and $\frac{1}{2}\sum_{I}(1+\varepsilon_{I}) = 2^{k} -1$ and the formula follows. If at least one of $d_{i}$ is odd, denote by $M$ the number
of $I$ with $d_{I}$ odd , then $\frac{1}{2}\sum_{I}(1+\varepsilon_{I}) = 2^{k} -1 - \frac{M}{2}$; we conclude by using $M = 2^{k-1}$.

The abelian varieties  $\prod_{I \in \mathcal{I}}\mathrm{Jac}(C_{I})$ and $\mathrm{Jac}(C)$ have the same dimension. Furthermore, if
$\eta_{j}= u^{j-1}du/v$ is a regular differential form on $C_{I}$ (for $1 \leq j \leq g_{I})$ then $\omega_{I,j}:= \phi^{\ast}_{I}(\eta_{j})=
x^{j-1}dx/y_{I}$ is regular on $C$  and these forms are linearly independent. Indeed an equality of type:
$$\sum_{I}\sum_{j=1}^{g_{I}}\lambda_{I,j}\omega_{I,j} = 0,$$
implies
$$\sum_{I}P_{I}(x)y_{I^{c}}=0,$$
where $I^{c}:= [1,k]\backslash I$ ; which implies the $P_{I}$'s are zero. The differential of $\Psi$ is an isomorphism, which proves that $\Psi$  is a separable isogeny.\\
\textbf{Remark.} We can deduce from the previous calculation, when $k \geq 2$, that the curve $C$ is not  hyperelliptic. In fact the canonical morphism  $P\mapsto (\omega_{I,j}(P))$ may be written $P \mapsto (1, x, x^2 , \ldots , y_{1}, \ldots , y_{k}, \ldots)$
thus is generically of degree 1, therefore an isomorphism (if $C$ were hyperelliptic, it would be a morphism of degree 2).
\begin{lem}Let $f_{1}, \ldots , f_{k}$ be polynomials  with coefficients in the finite field $\mathbb{F}_{q}$ such that the product $f_{1}... f_{k}$ is separable, and let $C$ be the associated curve. Let $A_{I} = q+1 - |C_{I}(\mathbb{F}_{q})|$ then,
$$|C(\mathbb{F}_{q})|= q +1 - \sum_{I \in \mathcal{I}}A_{I}.$$ \end{lem}
\textbf{Proof.} The abelian varieties $\mathrm{Jac}(C)$ and $\prod_{I}\mathrm{Jac}(C_{I})$ are $\mathbb{F}_{q}$-isogenous therefore have the same number of points on $\mathbb{F}_{q^{m}}$. If $$|C(\mathbb{F}_{q^{m}})| = q^{m}+ 1 - (\beta^{m}_{1}+ \ldots + \beta^{m}_{2g})$$ then $$|\mathrm{Jac}(C)(\mathbb{F}_{q^{m}})|= \prod_{1\leq i \leq 2g}(1-\beta^{m}_{i})$$ and if $$|C_{I}(\mathbb{F}_{q^{m}})| = q^{m}+ 1 - ((\alpha^{(I)}_{1})^{m}+ \ldots + (\alpha^{(I)}_{2g})^{m})$$ then $$|\mathrm{Jac}(C_{I})(\mathbb{F}_{q^{m}})|= \prod_{1\leq i \leq 2g}(1-(\alpha^{(I)}_{i})^{m})$$ and thus therefore
$$|\mathrm{Jac}(\mathcal{C})(\mathbb{F}_{q^{m}})|= \prod_{I}|\mathrm{Jac}(C_{I})(\mathbb{F}_{q^{m}})|= \prod_{I}\prod_{j}(1-(\alpha^{(I)}_{j})^{m}),$$
 $$|C_{I}(\mathbb{F}_{q^{m}})| = q^{m}+ 1 - ((\alpha^{(I)}_{1})^{m}+ \ldots + (\alpha^{(I)}_{2g})^{m}) = q^{m}+1- \sum_{I}\sum_{1\leq j\leq g_{I}}(\alpha_{j}^{(I)})^{m},$$
and, in particular , for $m = 1$, we get
$$|C(\mathbb{F}_{q})|= q +1 - \sum_{I}\sum_{1\leq j\leq g_{I}}\alpha_{j}^{(I)}= q + 1 - \sum_{I \in \mathcal{I}}A_{I} .$$
\newpage
\section{Numerical examples}
 Examples of genus-$g$ curves $y_1^2 = f_1$, $y_2^2 = f_2$
        over ${\bf F}_q$ having many points. The meaning of the quantities
         $A_1$, $A_2$, $A_3$, and $N$ is explained in the text. The following examples improve some results listed on the website  \texttt{http:// www.manypoints.org} and were computed using the computer algebra package Magma.\\

\begin{tabular}{llll}
\toprule
$q=17$ & $f_1 = x^4 + x^3 + 16 x^2 + 15 x + 1$   & $A_1 = -8$  & $N = 48 \in [\ldots,53]$ \cr
       & $f_2 = x^4 + 13 x^3 + 16 x^2 + 15$      & $A_2 = -6$  & \cr
       &                                         & $A_3 = -16$ & \cr
\midrule
$q=19$ & $f_1 = x^4 + x^3 + 18 x^2 + 13 x + 14$  & $A_1 = -8$  & $N = 52 \in [\ldots,60]$ \cr
       & $f_2 = x^4 + 4 x^3 + 18 x^2 + 7 x + 12$ & $A_2 = -8$  & \cr
       &                                         & $A_3 = -16$ & \cr
\midrule
$q=23$ & $f_1 = x^4 + 19 x^3 + 7$                & $A_1 = -9$  & $N = 60 \in [\ldots,67]$ \cr
       & $f_2 = x^3 + x + 11$                    & $A_2 = -9$  & \cr
       &                                         & $A_3 = -18$ & \cr
\midrule
$q=29$ & $f_1 = x^4 + x^3 + 28x^2+28x +18$                & $A_1 = -10$  & $N = 72 \in [\ldots,80]$ \cr
       & $f_2 = x^4 + 27x^3 + 27x^2+28x +15$                    & $A_2 = -10$  & \cr
       &                                         & $A_3 = -22$ & \cr
\midrule
$q=31$ & $f_1 = x^4 + x^3 + 30x^2+ 30x + 10$                & $A_1 = -11$  & $N = 76 \in [\ldots,84]$ \cr
       & $f_2 = x^4 + 30x^3 + 5x^2+ 19x +14$                    & $A_2 = -10$  & \cr
       &                                         & $A_3 = -23$ & \cr
\midrule
$q=37$ & $f_1 = x^4 + x^3 + 35x^2+ 32x +1$                & $A_1 = -12$  & $N = 88 \in [\ldots,96]$ \cr
       & $f_2 = x^4 + 11x^3 + 29x^2+ 13x +29$                    & $A_2 = -12$  & \cr
       &                                         & $A_3 = -26$ & \cr
\midrule
$q=41$ & $f_1 = x^4 + x^3 + 40x^2+ 40x + 36$                & $A_1 = -12$  & $N = 94 \in [\ldots,102]$ \cr
       & $f_2 = x^4 + 23x^3 + 10x^2+ 20x + 36$                    & $A_2 = -11$  & \cr
       &                                         & $A_3 = -29$ & \cr
\midrule
$q=43$ & $f_1 = x^4 + x^3 + 41x^2+ 34x + 1$                & $A_1 = -13$  & $N = 100 \in [\ldots,106]$ \cr
       & $f_2 = x^4 + 12x^3 + 17x^2+ 41x + 38$                    & $A_2 = -13$  & \cr
       &                                         & $A_3 = -30$ & \cr
\midrule
$q=47$ & $f_1 = x^4 + 25x^3 + x^2+ 2x + 31$      & $A_1 = -12$  & $N = 102 \in [\ldots,113]$ \cr
       & $f_2 = x^3 + x + 38$                    & $A_2 = -13$  & \cr
       &                                         & $A_3 = -29$ & \cr
\midrule
$q = 53$ & $f_{1}(x) = x^4 + x^3 + 52x^2 + 47x + 1 $     & $A_{1}  = -14$ &  $N = 120 \in [...,124]$ \cr
         &$f_{2}(x) = x^4 + 16x^3 + 36x^2 + 18x + 46 $   & $A_{2} = -14 $ & \cr
         &                                               & $A_{3} = - 38$ & \cr
\end{tabular}

\begin{center}
\begin{table}
\begin{tabular}{llll}
\toprule
 $q = 59$ & $f_{1}(x) = x^4 + x^3 + 54x^2 + 6x + 1 $   & $A_{1}  = -15$ &  $N = 124 \in [...,133]$ \cr
          & $f_{2}(x) = x^4 + 13x^3 + 22x^2 + 3x + 9 $  & $A_{2} = -15$ & \cr
          &                                             & $A_{3} = - 34$ &\cr
\midrule

$q = 61$ & $f_{1}(x) = x^4 + x^3 + 58x^2 + 18x + 1 $     & $A_{1}  = -15$& $N = 126 \in [...,137]$ \cr
         & $f_{2}(x) = x^4 + 27x^3 + 5x^2 + 47x + 28 $   & $A_{2} = -14$ & \cr
          &                                               &$A_{3} = - 35$ \cr
\midrule
$q = 67$ & $f_{1}(x) = x^4 + x^3 + 66x^2 + 57x + 1 $ & $A_{1}  = -16$& $N = 136 \in [...,148]$\cr
         & $f_{2}(x) = x^4 + 2x^3 + 49x^2 + 24x + 1 $ & $A_{2} = -16 $\cr
         &                                             &$A_{3} = - 36$ \cr


\midrule
For $q = 71$ & $f_{1}(x) = x^4 + x^3 + x^2 + 44x + 1 $ &  $A_{1}  = -16$ & $N = 144 \in [...,152]$ \cr
             & $f_{2}(x) = x^4 + 9x^3 + 8x^2 + 21x + 64 $ & $A_{2} = -16 $ \cr
            &                                             &$A_{3} = - 40$ \cr
\midrule
For $q = 73$ & $f_{1}(x) = x^4 + x^3 + 66x^2 + 57x + 1 $ & $A_{1}  = -17$& $N = 148 \in [...,156]$ \cr
             & $f_{2}(x) = x^4 + 2x^3 + 49x^2 + 24x + 1 $ & $A_{2} = -17 $\cr
             &                                            & $A_{3} = - 40$\cr
\midrule
For $q = 79$ & $f_{1}(x) = x^4 + x^3 + 3x^2 + 7x + 1 $     & $A_{1}  = -16$ & $N = 156 \in [...,165]$ \cr
             & $f_{2}(x) = x^4 + 4x^3 + 5x^2 + 24x + 68 $   & $A_{2} = -16$ \cr
             &                                              & $A_{3} = - 36$ \cr
\midrule
For $q = 83$ & $f_{1}(x) = x^4 + x^3 + x^2 + 5x + 1 $        &  $A_{1}  = -18$ &  $N = 162 \in [...,172]$ \cr
             & $f_{2}(x) = x^4 + 72x^3 + 54x^2 + 29x + 36 $  & $A_{2} = -16 $ \cr
             &                                               & $A_{3} = - 44$ \cr
\midrule
For $q = 89$ & $f_{1}(x) = x^4 + x^3 + 3x^2 + 7x + 1 $    & $A_{1}  = -16$ & $N = 168 \in [136,180]$\cr
             & $f_{2}(x) = x^4 + 4x^3 + 5x^2 + 24x + 68 $ & $A_{2} = -16 $ \cr
             &                                            &  $A_{3} = - 36$ \cr

\midrule
$q = 97$ & $f_{1}(x) = x^4 + 8x^3 + 3x^2 + 23x + 1 $ & $A_{1}  = -19$ &  $N = 180 \in [...,193]$ \cr
             & $f_{2}(x) = x^4 + 9x^3 + 63x^2 + 28x + 91 $ & $ A_{2} = -19 $ \cr
             &                                             & $A_{3} = - 44$ \cr

\midrule
$q = 5^2$ & $f_{1}(x) = x^4 + x^2 + rx $    & $A_{1}  = -9$ & $N = 68 \in [...,72]$\cr
             & $f_{2}(x) = x^4 -2rx^3 + rx + 2 $ & $A_{2} = -9 $ \cr
             &  with $r^2 + r + 1 = 0$             &  $A_{3} = - 24$ \cr
\midrule
$q = 13^2$ & $f_{1}(x) = (x-2)(x^{2}-2)  $ &  $A_{1}  = -26$ & $N = 295 \in [232,300]$ \cr
               & $f_{2}(x) = x^{4} -4x^{2}+1  $ & $A_{2} = -21$ \cr
               &                                &$A_{3} = -78$ \cr
\midrule
$q=17^{2}$ & $f_{1}(x) = (x+2)(x+10)(x^{2}+6)  $ &  $A_{1}  = -33$& $N = 454 \in [376,460]$ \cr
            & $f_{2}(x) = (x^{2}+10)(x^{2}+6x+3)$ & $A_{2} = -29$  \cr
            &                                     & $A_{3} = -102$ \cr
\midrule
\end{tabular}
\caption{Examples of genus-$5$.}
\end{table}
\end{center}

\begin{table}[h]
\begin{center}
\begin{tabular}{llll}
\toprule
$q = 23$ & $f_{1}(x) = x^{5} + 12x^4 + 19x^3 + x + 2 $ &  $A_{1}  = -13$ ~~~~ $N = 66 \in [60,78]$ \cr
         & $f_{2}(x) = x^{3} + 12x^2 + 18x + 4 $       & $A_{2} = -9$ \cr
         &                                             & $A_{3} = - 20$ \cr
\midrule
$q = 31$ & $f_{1}(x) = x^{5} + 6x^4 + 4x^3  + 4x^2 + 17x + 16 $ & $A_{1}  = -19$ ~~~~ $N = 84 \in [80,92]$ \cr
         & $f_{2}(x) = (x)^{3} + 13x^2 + 15x + 13 $ &    $A_{2} = -11$ \cr
         &                                          &     $A_{3} = -22$ \cr

\midrule
$q = 41$ & $f_{1}(x) = x^{5} + 31x^4 + 11x^3  + 14x^2 + 35x + 40 $ & $A_{1}  = - 23$~~~~ $N = 104 \in [102, 114]$ \cr
         & $f_{2}(x) = x^{3} + 23x^2 + 5x + 17 $ & $A_{2} = -12$ \cr
         &                                       &  $A_{3} = -27$ \cr

\midrule
$q = 59$ & $f_{1}(x) = x^{5} + 57x^4 + 17x^3  + 48 $ & $A_{1}  = -25$ ~~~~ $N = 134 \in [132, 150]$\cr
         & $f_{2}(x) = x^{3} + 2x + 22 $ & $A_{2} = -15 $ \cr
         &                               & $A_{3} = -34$ \cr

\bottomrule
\end{tabular}
\end{center}
\vskip1ex
\caption{Examples of genus-$6$ .}
\end{table}

\begin{table}[h]
\begin{center}
\begin{tabular}{llll}
\toprule
$q = 29$ & $f_{1}(x) = x^{6} + 9x^5 + 3x^4  + 1 $ & $A_{1}  = -46$ & $N = 80 \in [72, 100]$\cr
         & $f_{2}(x) = x^4 + 13x^{3} + 28x^2 + 12x +1 $ & $A_{2} = -40 $ \cr
         &                               & $A_{3} = -54$ \cr

\bottomrule
\end{tabular}
\end{center}
\vskip1ex
\caption{Examples of genus-$7$ .}
\end{table}


\subsection*{Acknowledgements}
We thank Professor Marc Hindry (Paris 7), the mathematics lab of the University Blaise Pascal in particular Nicolas Billerey and Marusia Rebelledo.\\

\end{document}